\documentclass[11pt,fullpage]{article}

\usepackage{amsmath}
\usepackage{amssymb}
\usepackage{amsfonts}
\usepackage{pstricks}
\usepackage{pst-plot}
\usepackage{xr}

\numberwithin{equation}{section}
\newtheorem{thm}{Theorem}[section] 
\newtheorem{prp}[thm]{Proposition}
\newtheorem{lmm}[thm]{Lemma}   
\newtheorem{crl}[thm]{Corollary}

\def\e_ref#1{(\ref{#1})}
\def\ov#1{\overline{#1}}

\def\lra{\longrightarrow}
\def\Lra{\Longrightarrow}

\def\io{\iota}
\def\ka{\kappa}
\def\om{\omega}

\def\De{\Delta}
\def\Si{\Sigma}

\def\PP{\Bbb{P}^2}
\def\PPP{\Bbb{P}^3}
\def\i{\infty}
\def\eset{\emptyset}

\begin{document}

\title{Completion of Katz-Qin-Ruan's Enumeration\\
of Genus-Two Plane Curves}
\author{Aleksey Zinger
\thanks{Partially supported by NSF Graduate Research Fellowship
and NSF grant DMS-9803166}}

\maketitle

\thispagestyle{empty}

\section{Introduction}

\noindent
In the past decade, significant progress has been 
in enumerative algebraic geometry based on ideas
of Gromov's compactness and quantum cohomology. 
In particular, \cite{KM} and~\cite{RT} derived
a recursive formula for the number~$N_d$ of
rational degree-$d$ plane curves passing through 
$(3d\!-\!1)$ points in general position.
In~\cite{I} and~\cite{P}, a simple relation between
the number~$N_{1,d}$ of fixed-$j$-invariant
elliptic degree-$d$ plane curves passing through
$(3d\!-\!1)$ points and the number $N_d$ is obtained.
The approaches in the two papers are drastically different.
In~\cite{P}, the number $N_{1,d}$ is computed by 
a beautiful degeneration argument.
In~\cite{I}, the number $N_{1,d}$ is compared
to the corresponding symplectic invariant as defined in~\cite{RT}.
Like the methods of \cite{KM} and~\cite{RT} in the genus-zero case,
the approach of \cite{I} applies to all projective spaces.\\

\noindent
The subject of this paper is the number~$N_{2,d}$ of 
genus-two degree-$d$ plane curves 
that have a fixed complex structure on the normalization and
pass through $(3d\!-\!2)$ points in general position.
Using a degeneration argument similar to~\cite{P},
\cite{KQR} express $N_{2,d}$
in terms of the numbers $N_{d'}$ \hbox{with $d'\!\le\! d$}.
Recently the author extended the approach of~\cite{I}
to obtain formulas for the genus-two numbers in $\PP$ and~$\PPP$.
However, the formulas for $N_{2,d}$ in~\cite{KQR} and~\cite{Z} 
are not equivalent.
The relation between the two is
$$N_{2,d}^Z=6\big(N_{2,d}^{KQR}+T_d\big),$$
where $T_d$ is the number of degree-$d$
tacnodal rational plane curves passing through $(3d\!-\!2)$ points.
The formulas in~\cite{Z} satisfy all the required classical checks 
that the author is aware~of.
In particular, $N_{2,4}^Z$ is the same as 
the corresponding number for three points and seven lines in~$\PPP$.
The author then explored the details of the argument~\cite{KQR}
and found three errors, one of which is significant.
They are described briefly in the paragraph following the table
and in more detail in Section~\ref{comp_sec}.
Once these errors are corrected, the formula of~\cite{Z} is recovered:

\begin{thm}
\label{main_thm}
$$N_{2,d}=3(d^2-1)N_d +
      \frac{1}{2} \sum_{d_1+d_2=d}\Big(d_1^2 d_2^2+28- 
    16\frac{9d_1d_2-1}{3d-2}\Big)
\binom{3d-2}{3d_1-1}d_1d_2N_{d_1}N_{d_2}.$$
\end{thm}

\noindent
The table below gives the numbers $N_{2,d}$ for small values of~$d$,
computed directly from Theorem~\ref{main_thm}.
The first three values have long been known to be zero.
We use $N_1\!=\!N_2\!=\!1$, $N_3\!=\!12$, $N_4\!=\!620$, $N_5\!=\!87,304$,
$N_6\!=\!26,312,976$, and $N_7\!=\!14,616,808,192$.

\begin{center}
\begin{tabular}{||c|c|c|c|c|c|c|c||}
\hline\hline
$d$&        1& 2& 3&    4&       5&        6& 7\\
\hline
$N_{2,d}$&  0& 0& 0& 14,400&  6,350,400& 3,931,128,000&
               3,718,909,209,600\\
\hline\hline
\end{tabular}
\end{center}

\vspace{.15in}

\noindent
The first step in the proof of Theorem~\ref{main_thm} via 
the recipe of~\cite{KQR} is Lemma~\ref{str_lmm},
which allows one to reduce the computation to 
a very degenerate genus-two curve.
The relevant intersection number is then
computed by Propositions~\ref{kqr_contr}-\ref{no_contr2}.
Propositions~\ref{kqr_contr} and~\ref{no_contr1}
are proved in~\cite{KQR}.
Proposition~\ref{no_contr2} is implied
by Remark~3.12 in~\cite{KQR}.
However, this remark is stated without a proof
and contradicts Propositions~\ref{extra_contr}.
This is the significant error in~\cite{KQR}.
A minor error is the statement about boundary relations
at the beginning of the proof of Lemma~2.18.
A~posteriori,  it turns out that this statement is in fact correct, 
at~least in the relevant cases,
but it does not follow from the argument given.
The remaining error is dividing by an extra factor of
six when computing contributions to the intersection number.\\

\noindent
Since our goal is to correct the computation
in~\cite{KQR}, we attempt to follow their notation 
as closely as possible.
The outline of this paper is as follows.
We first review the notation and setup in~\cite{KQR}.
In Section~\ref{comp_sec}, four propositions
are stated which imply~Theorem~\ref{main_thm}.
The last two sections prove the two propositions
not proved in~\cite{KQR}.
\\

\noindent
The author would like to thank T.~Mrowka for 
many discussions and encouragement.
He is also grateful to A.~J.~de~Jong,
J.~Starr, and R.~Vakil for their help with algebraic geometry.
In particular, it was A.~J.~de~Jong's idea 
to approach Corollary~\ref{contr_lmm2c} via the family of curves
of Lemma~\ref{contr_lmm2}.
Finally, the author thanks R.~Pandharipande for explaining
details of his argument in~\cite{P}
and Z.~Qin for careful consideration of the issues with~\cite{KQR}
raised by the author.

\section{Review of Notation and Setup}

\noindent
Denote by $\ov{\frak M}_2$ the Deligne-Mumford moduli space
of stable genus-two curves.
If  $d\!\ge\!3$, let
$$\ov{\frak M}_2(d)\equiv
\ov{\frak M}_{2,3d-2}\big(\PP,d\ell\big)$$
be Kontsevich's moduli space of stable maps
from $(3d\!-\!2)$-pointed genus-two curves to $\PP$ 
of degree~$d$, where $\ell\!\in\! H_2(\PP;\Bbb{Z})$
is the homology class of a line.
Let  $\pi\!:\ov{\frak M}_2(d)\!\lra\!\ov{\frak M}_2$
be the forgetful map.
Denote by $W_2(d)\!\subset\!\ov{\frak M}_2(d)$ the 
subset of stable maps with irreducible domains and
by $\ov{W}_2(d)$ the closure of $W_2(d)$ in~$\ov{\frak M}_2(d)$.\\

\noindent
Every element of $\ov{\frak M}_2(d)$ can be written as 
$\big[\mu\!:(D,p_1,\ldots,p_{3d-2})\big]$,
where $D$ is a prestable genus-two curve, 
$\mu\!:D\!\lra\!\PP$ is a (holomorphic) map,
and $p_1,\ldots,p_{3d-2}\!\in\! D$ are the marked points.
There are natural evaluation maps
$$e_i\!: \ov{\frak M}_2(d)\lra\PP,\quad
e_i\big(\big[\mu\!:(D,p_1,\ldots,p_{3d-2})\big]\big)
=\mu(p_i),\qquad i=1,\ldots,3d-2.$$
Let ${\cal L}_i=e_i^*\big({\cal O}_{\PP}(1)\big)$ and
$$Z=\big[\ov{W}_2(d)\big]\cap
c_1^2({\cal L}_1)\cap\ldots\cap c_1^2({\cal L}_{3d-2})
\in H_6\big(\ov{W}_2(d)\big).$$
If $q_1,\ldots,q_{3d-2}$ are points in $\PP$ in general position,
then 
$\big\{e_1\!\times\!\ldots\!\times\!e_{3d-2}\big\}^{-1}
(q_1\!\times\!\ldots\!\times\! q_{3d-2})$
is a representative for~$Z$; see~\cite{KQR} for details.

\begin{lmm}
\label{str_lmm}
For every $[C]\!\in\!\ov{\frak M}_2$,
$$N_{2,d}=[\pi^{-1}(C)]\cdot Z,$$
where $[\pi^{-1}(C)]\cdot Z$ is 
the intersection pairing of $\pi^{-1}([C])$ and  $Z$ in
$\ov{W}_2(d)$.
\end{lmm}

\noindent
This is a special case of Lemma~2.5 in \cite{KQR}.
In particular, if $C_0$ consists of two rational components
identified at $3$ pairs of points, i.e.

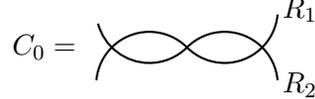
\begin{pspicture}(0,-1.7)(10,0)
\psarc(7.1,-1){.707}{45}{175}\psarc(7.1,0){.707}{195}{315}
\psarc(8.1,-1){.707}{5}{135}\psarc(8.1,0){.707}{225}{355}
\rput(5.7,-.5){$C_0=$}\rput(9.1,0){$R_1$}\rput(9.1,-1){$R_2$}
\end{pspicture}

\noindent
then $N_{2,d}\!=\![\pi^{-1}(C_0)]\cdot Z$.
The space $\pi^{-1}(C_0)\!\subset\!\ov{\frak M}_2(d)$ can be written 
as the disjoint union~$\bigsqcup W_T$,
where $W_T$ is the space of stable maps 
$\big[\mu\!:(D,p_1,\ldots,p_{3d-2})\big]$,
such that the domain $D$ is the union of $R_1$, $R_2$,
and trees $T_1,\ldots,T_s$ of $\Bbb{P}^1$ in a way encoded by~$T$.
The stable reduction of $D$ must be~$C_0$.
See Figure~1 below for some examples.\\

\noindent
In order to compute  $[\pi^{-1}(C_0)]\cdot Z$, 
\cite{KQR} consider the intersection of $Z$ with 
every nonempty space~$W_T$.
It is fairly easy to show that $Z\cap W_T$ is empty for all
but a small number of trees~$T$, independent of~$d$.
If $\big[\mu\!:(D,p_1,\ldots,p_{3d-2})\big]\!\in\! Z\cap W_T$,
the map $\mu\!:D\!\lra\!\PP$ has degree~$d$ and passes
through $3d\!-\!2$ points in $\PP$ in general position.
Thus, if $D_1,\ldots,D_m$ are the irreducible components of $D$
to which $\mu$ restricts non-trivially,
$m\!=\!1$ or $m\!=\!2$.
Then $D$ can have at most two components, other than $R_1$, $R_2$,
on which the map $\mu$ is constant.\\

\noindent
The complete list of possibilities for $D$, up to 
equivalence, is given in Figure~1.
Denote by $C_{ij}$ the curve as in the $i$th row
and $j$th column of Figure~1.
Similarly, denote by $W_{ij}$ be the space of stable maps  
with domain~$C_{ij}$ and a distribution of the degree~$d$
between the components of $C_{ij}$ such that
the image of some stable map in $W_{ij}$
passes through $(3d\!-\!2)$ points.
We clarify this statement in the relevant cases:\\
(1) if $\big[\mu\!:(D,p_1,\ldots,p_{3d-2})\big]$ lies in
$W_{13}$, $W_{32}$, $W_{41}$, $W_{43}$, or $W_{5j}$,
the degree of $\mu|D_i$ is $d_i\!\neq\!0$, and 
the restriction of $\mu$ to all other components is constant;\\
(2) if $\big[\mu\!:(D,p_1,\ldots,p_{3d-2})\big]$ lies in 
$W_{24}$, $W_{31}$, or $W_{42}$, 
the degree of $\mu|D_1$  is $d_1\!\neq\!0$,
$\mu|R_i$ is constant, and
in the case of $W_{42}$ the restriction of $\mu$
to the vertical component (in the diagram) is constant.\\
Furthermore, for stability reasons, 
every component of $C_{ij}$, on which $\mu$ is constant
and which does not contain three singular point of~$C_{ij}$,
must contain one of the marked points~$p_i$.

\begin{figure}
\begin{pspicture}(-1.1,-10)(10,0)
\psarc(.5,-1){.707}{45}{175}\psarc(.5,0){.707}{195}{315}
\psarc(1.5,-1){.707}{5}{135}\psarc(1.5,0){.707}{225}{355}
\rput(2.5,0){$R_1$}\rput(2.5,-1){$R_2$}
\psarc(4.5,-1){.707}{45}{175}\psarc(4.5,0){.707}{195}{315}
\psarc(5.5,-1){.707}{5}{135}\psarc(5.5,0){.707}{225}{355}
\psline(4.3,-.5)(5,.2)
\rput(6.5,0){$R_1$}\rput(6.5,-1){$R_2$}\rput(5.25,.15){$D_1$}
\psarc(8,-1){.707}{45}{110}\psarc(8,0){.707}{250}{315}
\psarc(9,-1){.707}{5}{135}\psarc(9,0){.707}{225}{355}
\psline(8,-1.2)(8,.2)
\rput(10,0){$R_1$}\rput(10,-1){$R_2$}\rput(8.35,.15){$D_1$}
\psarc(12.5,-1){.707}{45}{175}\psarc(12.5,0){.707}{195}{315}
\psarc(13.5,-1){.707}{5}{135}\psarc(13.5,0){.707}{225}{355}
\psline(12.3,-.5)(13,.2)\psline(12.3,.3)(13,-.2)
\rput(14.5,0){$R_1$}\rput(14.5,-1){$R_2$}\rput(12.4,.5){$D_1$}
%end of first row of pictures 
\psarc(0.5,-3){.707}{45}{175}\psarc(0.5,-2){.707}{195}{315}
\psarc(1.5,-3){.707}{5}{135}\psarc(1.5,-2){.707}{225}{355}
\psline(.6,-2.5)(1.3,-1.8)\psline(.4,-2.5)(-.3,-1.8)
\rput(2.5,-2){$R_1$}\rput(2.5,-3){$R_2$}
\rput(1.55,-1.85){$D_2$}\rput(.05,-1.7){$D_1$}
\psarc(4.5,-3){.707}{45}{175}\psarc(4.5,-2){.707}{195}{315}
\psarc(5.5,-3){.707}{5}{135}\psarc(5.5,-2){.707}{225}{355}
\psline(4.6,-2.5)(5.3,-1.8)\psline(4.3,-2.5)(3.8,-3.2)
\rput(6.5,-2){$R_1$}\rput(6.5,-3){$R_2$}
\rput(5.55,-1.85){$D_1$}\rput(4.2,-3.1){$D_2$}
\psarc(8,-3){.707}{45}{110}\psarc(8,-2){.707}{250}{315}
\psarc(9,-3){.707}{5}{135}\psarc(9,-2){.707}{225}{355}
\psline(7.9,-3.2)(7.9,-1.8)\psline(8.1,-2.5)(8.8,-1.8)
\rput(10,-2){$R_1$}\rput(10,-3){$R_2$}\rput(9.05,-1.85){$D_1$}
\psarc(12.5,-3){.707}{45}{110}\psarc(12.5,-2){.707}{250}{315}
\psarc(13.5,-3){.707}{5}{135}\psarc(13.5,-2){.707}{225}{355}
\psline(12.4,-3.2)(12.4,-1.8)\psline(12.2,-2.1)(13.1,-1.7)
\rput(14.5,-2){$R_1$}\rput(14.5,-3){$R_2$}\rput(13.4,-1.7){$D_1$}
%end of second row of pictures 
\psarc(.5,-5){.707}{45}{110}\psarc(.5,-4){.707}{250}{315}
\psarc(1.5,-5){.707}{5}{135}\psarc(1.5,-4){.707}{225}{355}
\psline(0.1,-4.6)(1,-3.9)\psline(0.1,-4.4)(1,-5.1)
\rput(2.5,-4){$R_1$}\rput(2.5,-5){$R_2$}\rput(1.27,-3.9){$D_1$}
\psarc(4.5,-5){.707}{45}{115}\psarc(4.5,-4){.707}{245}{315}
\psarc(5.5,-5){.707}{70}{135}\psarc(5.5,-4){.707}{225}{290}
\psline(4.4,-5.2)(4.4,-3.8)\psline(5.6,-5.2)(5.6,-3.8)
\rput(4.1,-3.8){$D_1$}\rput(4.8,-4.1){$R_1$}\rput(4.8,-4.9){$R_2$}
\psarc(8,-5){.707}{45}{175}\psarc(8,-4){.707}{195}{315}
\psarc(9,-5){.707}{5}{135}\psarc(9,-4){.707}{225}{355}
\psline(8,-4.5)(8,-3.6)\psline(7.8,-3.7)(8.6,-3.7)
\psline(7.4,-3.9)(8.2,-3.9)
\rput(10,-4){$R_1$}\rput(10,-5){$R_2$}
\rput(7.1,-3.9){$D_1$}\rput(8.9,-3.7){$D_2$}
\psarc(12.5,-5){.707}{45}{175}\psarc(12.5,-4){.707}{195}{315}
\psarc(13.5,-5){.707}{5}{135}\psarc(13.5,-4){.707}{225}{355}
\psline(12.6,-4.5)(13.3,-3.8)
\psline(12,-4.1)(13.1,-4.1)\psline(12.35,-4.25)(13.05,-3.55)
\rput(14.5,-4){$R_1$}\rput(14.5,-5){$R_2$}
\rput(12.25,-3.9){$D_1$}\rput(13.05,-3.4){$D_2$}
%end of third row of pictures 
\psarc(.5,-7){.707}{45}{110}\psarc(.5,-6){.707}{250}{315}
\psarc(1.5,-7){.707}{5}{135}\psarc(1.5,-6){.707}{225}{355}
\psline(.4,-7.2)(.4,-5.8)
\psline(.2,-6.1)(1.1,-5.7)\psline(.6,-6.9)(-.3,-7.3)
\rput(2.5,-6){$R_1$}\rput(2.5,-7){$R_2$}
\rput(1.35,-5.7){$D_1$}\rput(-.6,-7.3){$D_2$}
\psarc(4.5,-7){.707}{45}{110}\psarc(4.5,-6){.707}{250}{315}
\psarc(5.5,-7){.707}{5}{135}\psarc(5.5,-6){.707}{225}{355}
\psline(4.4,-7.2)(4.4,-5.8)
\psline(4.2,-6.1)(5.1,-5.7)\psline(4.55,-5.85)(5.6,-5.85)
\rput(6.5,-6){$R_1$}\rput(6.5,-7){$R_2$}
\rput(5.6,-5.62){$D_1$}
\psarc(8.5,-7){.707}{45}{110}\psarc(8.5,-6){.707}{250}{315}
\psarc(9.5,-7){.707}{60}{135}\psarc(9.5,-6){.707}{225}{305}
\psline(8.4,-7.2)(8.4,-5.8)\psline(9.6,-7.2)(9.6,-5.8)
\psline(8.5,-6)(7.7,-5.7)
\rput(7.45,-5.7){$D_1$}
\rput(8.8,-6.1){$R_1$}\rput(8.8,-6.9){$R_2$}
\psarc(12.5,-7){.707}{45}{110}\psarc(12.5,-6){.707}{250}{315}
\psarc(13.5,-7){.707}{5}{135}\psarc(13.5,-6){.707}{225}{355}
\psline(12.4,-7.2)(12.4,-5.8)
\psline(13.6,-6.5)(14.2,-5.7)
\psline(13,-6)(14.05,-6)
\rput(14.5,-6.1){$R_1$}\rput(14.5,-7){$R_2$}
\rput(13.4,-5.79){$D_2$}\rput(14.5,-5.6){$D_1$}
%end of fourth row of pictures 
\psarc(.5,-9){.707}{45}{110}\psarc(.5,-8){.707}{250}{315}
\psarc(1.5,-9){.707}{5}{135}\psarc(1.5,-8){.707}{225}{355}
\psline(.4,-9.2)(.4,-7.8)\psline(.5,-8.5)(-.6,-8.7)
\psline(0,-8.75)(0,-7.8)\psline(-.4,-8.5)(-.45,-9.4)
\rput(2.5,-8){$R_1$}\rput(2.5,-9){$R_2$}
\rput(-.3,-7.85){$D_2$}\rput(-.15,-9.35){$D_1$}
\psarc(4.5,-9){.707}{45}{110}\psarc(4.5,-8){.707}{250}{315}
\psarc(5.5,-9){.707}{5}{135}\psarc(5.5,-8){.707}{225}{355}
\psline(4.4,-9.2)(4.4,-7.8)\psline(4.5,-8.5)(3.7,-8.65)
\psline(4,-8.75)(4,-7.8)\psline(4.15,-8.2)(3.05,-8.4)
\rput(6.5,-8){$R_1$}\rput(6.5,-9){$R_2$}
\rput(3.4,-8){$D_2$}\rput(4,-7.6){$D_1$}
\psarc(8.5,-9){.707}{45}{110}\psarc(8.5,-8){.707}{250}{315}
\psarc(9.5,-9){.707}{60}{135}\psarc(9.5,-8){.707}{225}{305}
\psline(8.4,-9.2)(8.4,-7.8)\psline(9.6,-9.2)(9.6,-7.8)
\psline(8.5,-8)(7.7,-7.7)\psline(9.5,-8)(10.3,-7.7)
\rput(7.45,-7.7){$D_1$}\rput(10.55,-7.7){$D_2$}
\rput(8.8,-8.1){$R_1$}\rput(8.8,-8.9){$R_2$}
\psarc(12.5,-9){.707}{45}{110}\psarc(12.5,-8){.707}{250}{315}
\psarc(13.5,-9){.707}{60}{135}\psarc(13.5,-8){.707}{225}{305}
\psline(12.4,-9.2)(12.4,-7.8)\psline(13.6,-9.2)(13.6,-7.8)
\psline(12.5,-8)(11.7,-7.7)\psline(12.3,-7.75)(11.9,-8.5)
\rput(11.45,-7.7){$D_1$}
\rput(12.8,-8.1){$R_1$}\rput(12.8,-8.9){$R_2$}
\rput(11.75,-8.65){$D_2$}
\rput(6.6,-10){Figure~1}
\end{pspicture}
\end{figure}

\section{Computation of the Intersection Number}
\label{comp_sec}

\begin{prp}
\label{kqr_contr}
The contribution to $[\pi^{-1}(C_0)]\cdot Z$ from $W_{11}$ is 
$$\frac{3(d-1)(d-2)(d-3)}{d}N_d+
\frac{1}{2}\sum_{d_1+d_2=d}
\Big( d_1^2 d_2^2 -6d_1d_2 -4 +18 \frac{d_1d_2}{d}\Big)
\binom{3d-2}{3d_1-1}d_1d_2N_{d_1}N_{d_2}.$$
\end{prp}

\noindent
This proposition is essentially proved in~\cite{KQR};
see equation~(2.9) and Lemmas~2.12, 2.16, and~3.2 in~\cite{KQR}.
The above number is six times the number given
by Theorem~1.1 of~\cite{KQR}.
It is easy to see that the authors divide by six an extra time.
For example, in Lemma~2.12, one should take {\it ordered}
triplets of nodes, 
i.e.~$\binom{d_1d_2}{3}$ should be replaced~by
$$d_1d_2(d_1d_2-1)(d_1d_2-2),$$ 
since they are dividing  by the order of~$\hbox{Aut}(C_0)$.
Similarly, the number in Lemma~2.16 should be replaced by six times
itself.

\begin{prp}
\label{extra_contr}
The contribution to $[\pi^{-1}(C_0)]\cdot Z$ from $W_{13}$ is 
$$\frac{6(3d^2-12d+9)n_d}{d} + 
3\sum_{d_1+d_2=d} 
\Big(d_1d_2+4-9\frac{d_1d_2}{d}\Big)
\binom{3d-2}{3d_1-1}d_1d_2N_{d_1}N_{d_2}.$$
\end{prp}

\noindent
We prove this proposition in Section~\ref{extra_contr_sec}.
What we show is that $\ov{W}_2(d)\cap W_{13}$
is the space of all stable maps 
$\big[\mu\!:(D,p_1,\ldots,p_{3d-2})\big]$
such that $\mu(D)$ is a tacnodal curve in $\PP$,
and $\mu$ maps the two nodes of $D$ to the same tacnode 
of~$\mu(D)$.
The number of Proposition~\ref{extra_contr} is~$6T_d$.
Note that the number~$T_d$ is well-known; see equation~(1.2) in~\cite{DH}
and Subsection~3.2 in~\cite{V1}.

\begin{prp}
\label{no_contr1}
If $(i,j)\!\in\!\big\{(1,2),(1,4),(2,1),(2,2),(2,3),(3,3),(3,4),(4,4)\big\}$,
$Z\cap W_{ij}\!=\!\eset$.
Thus, $W_{ij}$ does not contribute to~$[\pi^{-1}(C_0)]\cdot Z$.
\end{prp}

\noindent
Most of this proposition is proved
by  Lemmas~2.18 and~3.7 of~\cite{KQR}.
The cases $(i,j)\!=\!(3,3)$ and $(i,j)\!=\!(3,4)$
can be deduced from the proofs of these two lemmas.
The modification required is similar
to the extension of the main part of the proof of Lemma~1 in~\cite{P}
to cases of multiple blowups;
see also the proof of Lemma~\ref{nc_lmm3} below.
Note that since Lemma~3.7 of~\cite{KQR} does not apply to the remaining
possibilities for~$(i,j)$, neither does Lemma~2.18 of~\cite{KQR}.

\begin{prp}
\label{no_contr2}
If $(i,j)\!\in\!\big\{(2,4),(3,1),(3,2),(4,1),(4,2),(4,3),(5,1),(5,2),
(5,3),(5,4)\big\}$,\\ \hbox{$Z\cap W_{ij}\!=\!\eset$}.
Thus, $W_{ij}$ does not contribute to~$[\pi^{-1}(C_0)]\cdot Z$.
\end{prp}

\noindent
We prove this proposition in the next section.
The number in Theorem~\ref{main_thm} is the sum of the numbers
in Propositions~\ref{kqr_contr} and~\ref{extra_contr}.
However, one has to make use of Kontsevich's recursion
to obtain the formula in Theorem~\ref{main_thm}:
$$N_d=\frac{1}{6(d-1)} 
\sum_{d_1+d_2=d} \Big(d_1d_2-2\frac{(d_1-d_2)^2}{3d-2}\Big)
\binom{3d-2}{3d_1-1}d_1d_2N_{d_1}N_{d_2}.$$

\section{Proof of Proposition~\ref{no_contr2}}

\subsection{The Semi-Standard Cases}

\noindent
We prove Proposition~\ref{no_contr2} by exhibiting
conditions that stable maps in $\ov{W}_2(d)\cap W_{ij}$ must satisfy.
This approach is analogous to methods in~\cite{P}
and~\cite{KQR}, but we make no use of the spaces
$X$ and~$Y$ of these two papers.
It should be possible to describe  the space 
\hbox{$\ov{W}_2(d)\cap\pi^{-1}(C_0)\!\subset\!\ov{\frak M}_2(d)$}
explicitly by using arguments as in this section 
to obtain necessary conditions for an  element of $\pi^{-1}(C_0)$
to be in $\ov{W}_2(d)$
and by applying methods similar to the next section
to show that these conditions are sufficient.
However, much less is needed to prove Theorem~\ref{main_thm}.
\\

\noindent
Suppose $\big[\mu\!:(D,p_1,\ldots,p_{3d-2})\big]
\!\in\!\ov{W}_2(g)\cap W_{ij}$.
Then by definition of stable-map convergence, 
there exist\\
(T1) a one-parameter family of curves
$\tilde{\eta}\!:\tilde{\cal F}\!\lra\!\De$
such that $\De$ is a neighborhood of~$0$ in~$\Bbb{C}$,
$\tilde{\cal F}$ is a smooth space, $\tilde{\eta}^{-1}(0)\!=\!D$,
and $C_t\!\equiv\!\eta^{-1}(t)$ is a smooth genus-two curve
for all $t\!\in\!\De^*\!\equiv\!\De\!-\!{0}$;\\
(T2) a map $\tilde{\mu}\!:\tilde{\cal F}\!\lra\!\PP$ such that
$\tilde{\mu}|\eta^{-1}(0)\!=\!\mu$.\\
In many cases, $\tilde{\cal F}$ can be obtained 
by a sequence of blowups from another smooth bundle
$\eta\!:{\cal F}\!\lra\!\De$ of curves.
This observation is used often in the proofs of the lemmas that follow.

\begin{lmm}
\label{nc_lmm1}
If $\big[\mu\!:(D,p_1,\ldots,p_{3d-2})\big]\!\in\!\ov{W}_2(d)\cap W_{24}$
and the degree of $\mu|D_1$ is~$d$,
$\mu(D)$ has a cusp at~$\mu(p_i)$ for some $i\!=\!1,\ldots,3d\!-\!2$.
\end{lmm}

\noindent
{\it Proof:}
(1) Let $\tilde{\eta}\!:\tilde{\cal F}\!\lra\!\De^*$
be a family as in (T1) above with central fiber $\tilde{C}_0\!=\!D$,
and $\tilde{\mu}\!:\tilde{\cal F}\!\lra\!\PP$ a map as in (T2).
Then there exists another family
$\eta\!:{\cal F}\!\lra\!\De$ as in (T1)
such that the central fiber is $C_{13}$ 
and $\tilde{\cal F}$ is the blowup of ${\cal F}$
at a smooth point $p\!\in\!D_1\subset C_{13}$.\\
(2) Let $\psi\!\in\!H^0\big(C_{13};\om_{C_{13}}\big)$
be an element such that $\psi|D_1\!\neq\!0$. 
From the point of view of complex geometry, 
$H^0\big(C_{13};\om_{C_{13}}\big)$ is the space harmonic 
$(1,0)$-forms on the three components of~$C_{13}$,
which have simple poles at the singular points
with residues that add up to zero at each node.
Thus, such an element exists.
Let $(t,w)$ be coordinates near $p\!\in\!{\cal F}$
such that $w$ is the vertical coordinate, 
i.e.~$d\eta\big|\frac{\partial}{\partial w}\!=\!0$.
Then $\psi$ extends to a family of elements 
$\psi_t\!\in\!H^0\big(C_t;\om_{C_t}\big)$ such that
\begin{equation}\label{nc_lmm1_e2}
\psi_t\big|_w=a\big(1+o(1_{(t,w)})\big)dw,
\end{equation}
for some $a\!\in\!\Bbb{C}^*$.\\
(3) On a neighborhood of~$D_1^*\!\subset\! D\!=\!C_{24}$, 
the complement of the node in~$D_1$,
we have local coordinates $(t,z)\!\lra\!\big(t,w\!=\!tz,[1,z]\big)$.
Note that in these coordinates, \e_ref{nc_lmm1_e2} becomes
\begin{equation}\label{nc_lmm1_e2b}
\psi_t\big|_z=at\big(1+o(1_t)\big)dz.
\end{equation}
Let $L_1$ and $L_2$ be any two lines in general position in~$\PP$.
In particular, we assume that they miss the image under $\mu$
of the node of~$D_1$.
Then for all $t\!\in\!\De$, sufficiently small,
\begin{equation}\label{nc_lmm1_e3}
\mu_t^{-1}(L_i)=\big\{z^{(i)}_1(t),\ldots,z^{(i)}_d(t)\big\}
\subset C_t \quad\hbox{and}\quad
z^{(i)}_j(t)=z^{(i)}_j(0)+o(1_t),
\end{equation}
where $\mu_t\!=\!\tilde{\mu}|C_t$.
Since $\sum z^{(1)}_j(t)$  and $\sum z^{(2)}_j(t)$
are linearly equivalent divisors in~$C_t$,
\begin{equation}\label{nc_lmm1_e5}
\sum_{j=1}^{j=d}\int_{z^{(1)}_j(t)}^{z^{(2)}_j(t)}\psi_t=0
\qquad\forall t\in\De^*,
\end{equation}
where the line integrals are taken inside of the coordinate chart.
Plugging \e_ref{nc_lmm1_e2b} and \e_ref{nc_lmm1_e3}
into \e_ref{nc_lmm1_e5} gives
\begin{equation}\label{nc_lmm1_e7}
at\sum_{j=1}^{j=d}\big(z^{(2)}_j(0)-z^{(1)}_j(0)+o(1_t)\big)=0
\qquad\forall t\!\in\!\De^*.
\end{equation}
Dividing this equation by $at$ and then taking the limit as $t\!\lra\!0$,
we conclude that
\begin{equation}\label{nc_lmm1_e9}
\sum_{j=1}^{j=d}z^{(1)}_j(0)=\sum_{j=1}^{j=d}z^{(2)}_j(0).
\end{equation}
Condition  \e_ref{nc_lmm1_e9} can be explicitly interpreted as follows.
Let $[u,v]$ be homogeneous coordinates on $D_1$
such that $z\!=\!\frac{v}{u}$.
Then a map  $D_1\!\lra\!\PP$ of degree-$d$
corresponds to three homogeneous polynomials
$$p_i=\sum_{j=0}^{j=d}p_{ij}u^jv^{d-j}.$$
Since equality \e_ref{nc_lmm1_e9} holds for a dense subset of lines
in $\PP$, there exists $K\!=\!K(\mu)\!\in\!\Bbb{P}^1$ such~that
\begin{gather}
\frac{c_0p_{0,d-1}+c_1p_{1,d-1}+c_2p_{2,d-1}}
{c_0p_{0,d}+c_1p_{1,d}+c_2p_{2,d}}
=K
\quad\forall (c_0,c_1,c_2)\!\in\!\Bbb{C}^3\!-\!\{0\}\Lra\notag\\
\label{nc_lmm1_e11}
(p_{0,d-1},p_{1,d-1},p_{2,d-1})=
K(p_{0,d},p_{1,d},p_{2,d}).
\end{gather}
Equation~\e_ref{nc_lmm1_e11} imposes two linearly independent
conditions on the map~$\mu|D_1$
if $\mu\!\in\!\ov{W}_2(d)\cap W_{24}$.
Geometrically, they mean that $\mu\big(D\big)$
has a cusp at the image of the node of~$D_1$.

\begin{crl}
\label{nc_lmm1c}
If $\big[\mu\!:(D,p_1,\ldots,p_{3d-2})\big]\!\in\!Z\cap W_{24}$,
the degree of $\mu|D_1$ is less than~$d$.
\end{crl}

\noindent
{\it Proof:} Suppose the degree of $\mu|D_1$ is $d$.
Then by Lemma~\ref{nc_lmm1}, $\mu(D_1)$ has a cusp at 
the image of the node of~$D_1$.
Since the points $q_1,\ldots,q_{3d-2}$ are in general position,
$\mu(D_1)$ has one simple cusp and $\binom{d-1}{2}\!-\!1$ simple nodes.
Let $\tilde{\cal F}$ and $\tilde{\mu}$ be as in 
the proof of Lemma~\ref{nc_lmm1}.
Then $\tilde{\mu}(C_t)$ converges to~$\mu(D_1)$.
By Lemma~2.4.1 or Example~3.2.2 in~\cite{V2}, 
$D_1$ must have an elliptic tail,
i.e.~the map $\tilde{\mu}\!:\tilde{\cal F}\!\lra\!\PP$
cannot exist.
In the given case, this can also be seen directly as follows.
The image under $\mu_t$ of the intersection of $C_t$ with the coordinate
chart described in (3) of the proof of Lemma~\ref{nc_lmm1}
has $\binom{d-1}{2}\!-\!1$ simple nodes, 
close to the simple nodes of~$\mu(D_1)$.
The complement of the coordinate chart in $C_t$ 
is a genus two curve with a small coordinate neighborhood removed.
Thus, it contributes at least $2$ to the arithmetic genus of~$\mu(C_t)$.
This means that the arithmetic genus of $\mu(C_t)$
is at least  $\binom{d-1}{2}\!+\!1$, instead of~$\binom{d-1}{2}$.

\begin{lmm}
\label{nc_lmm2}
The image of every element 
$\big[\mu\!:(D,p_1,\ldots,p_{3d-2})\big]\!\in\!\ov{W}_2(d)\cap W_{43}$
has a cusp at~$\mu(p_i)$ for some $i\!=\!1,\ldots,3d\!-\!2$.
The same is true  for every element of $\ov{W}_2(d)\cap W_{42}$
such that the degree of $\mu|D_1$ is~$d$.
Thus, $Z\cap W_{43}\!=\!\eset$, while
for every element 
$\big[\mu\!:(D,p_1,\ldots,p_{3d-2})\big]\!\in\! Z\cap W_{42}$,
the degree of $\mu|D_1$ is less than~$d$.
\end{lmm}

\noindent
{\it Proof:}
(1) The proof of the first statement is nearly the same 
as the proof of Lemma~\ref{nc_lmm1}.
The only difference is that the central fiber of ${\cal F}$
will be~$C_{32}$.\\
(2) The family $\tilde{\cal F}$ of the second claim of this
lemma is obtained from $\tilde{\cal F}$ of Lemma~\ref{nc_lmm1}
by blowing up a smooth point of the exceptional 
divisor $D_1\!\subset\! C_{24}$.
Thus, nearly the same argument as in  Lemma~\ref{nc_lmm1}
applies if the degree of $\mu|D_1$ is~$d$;
see~\cite{P} for an extension in an analogous situation.

\begin{lmm}
\label{nc_lmm3}
If $(i,j)\!\in\!\big\{(5,2),(5,4)\big\}$,
the image of every element of 
$\ov{W}_2(d)\cap W_{ij}$ is a two-component rational cuspidal curve. 
The same is true for all
$\big[\mu\!:(D,p_1,\ldots,p_{3d-2})\big]\!\in\! \ov{W}_2(d)\cap W_{42}$,
such that the degree of $\mu|D_1$ is less than~$d$.
Thus, $Z\cap W_{ij}\!=\!\eset$ in all three cases.
\end{lmm}

\noindent
{\it Proof:}
(1) We first consider the case
 $\big[\mu\!:(D,p_1,\ldots,p_{3d-2})\big]\!\in\! \ov{W}_2(d)\cap W_{42}$
and the degree of $\mu|D_1$ is $d_1\!<\!d$.
The case $d_1\!=\!d$ is considered in Lemma~\ref{nc_lmm2}.
The family $\tilde{\cal F}\!\lra\!\De$ corresponding to this
case can be obtained as follows.
We start with a family ${\cal F}\!\lra\!\De$ as in (2)
of the proof of Lemma~\ref{nc_lmm1},
blow it up at a smooth point  $p\!\in\! D_1\!\subset\! C_{13}$,
and then blow up the resulting space
at a smooth point $p_1$ of the new exceptional divisor 
$E\!\equiv\! D_1\!\subset\! C_{24}$.
Denote the last exceptional divisor by~$E_1$.
We use coordinates $(t,z)$ near $E^*$ as before
and coordinates 
$(t,z_1)\!\lra\!\big(t,z=p_1\!+\!tz_1,[1,z_1]\big)$
near $E_1^*$.
Then,
\begin{gather*}
\psi_t\big|_z=at\big(1+o(1_t)\big)dz,\qquad
\psi_t\big|_{z_1}=at^2\big(1+o(1_t)\big)dz_1;\\
\mu_t^{-1}(L_i)=\big\{ z_{1,1}^{(i)}(t),\ldots,z_{1,d_1}^{(i)}(t),
z_{d_1+1}^{(i)}(t),\ldots,z_d^{(i)}(t)\big\}\subset C_t,
\qquad\hbox{with}\\
z_{1,j}^{(i)}(t)=z_{1,j}^{(i)}(0)+o(1_t),\quad
z_j^{(i)}(t)=z_j^{(i)}(0)+o(1_t);\\
\sum_{j=1}^{j=d_1}\int_{z^{(1)}_{1,j}(t)}^{z^{(2)}_{1,j}(t)}\psi_t+
\sum_{j=d_1+1}^{j=d}\int_{z^{(1)}_j(t)}^{z^{(2)}_j(t)}\psi_t=0
\qquad\forall t\!\in\!\De^*.
\end{gather*}
Each line integral is taken inside the corresponding coordinate chart.
Proceeding as in the proof of Lemma~\ref{nc_lmm1}, we obtain
\begin{gather*}
at^2\sum_{j=1}^{j=d_1}\big(z^{(1)}_{1,j}(0)-z^{(2)}_{1,j}(0)+o(1_t)\big)
+at \sum_{j=d_1+1}^{j=d}\big(z^{(2)}_j(0)-z^{(1)}_j(0)+o(1_t)\big)=0
\qquad\forall t\!\in\!\De^*\\
\Lra\sum_{j=d_1+1}^{j=d}z^{(1)}_j(0)=
\sum_{j=d_1+1}^{j=d}z^{(2)}_j(0).
\end{gather*}
As before, the last identity implies that $\mu|E$
maps $z\!=\!\i\!\in\! E$ to a cusp of~$\mu(E)$.\\
(2)  The argument in the case of $W_{54}$ is the same,
except we replace the family ${\cal F}$
of Lemma~\ref{nc_lmm1} with the family ${\cal F}$ 
of (1) of Lemma~\ref{nc_lmm2}.
Finally, the case of $W_{52}$
simply involves an extra blowup at a smooth point
as compared to the case of~$W_{42}$.

\begin{lmm}
\label{nc_lmm4}
If $(i,j)\!\in\!\big\{(4,1),(5,1)\big\}$,
the image of every element of 
$\ov{W}_2(d)\cap W_{ij}$ is a two-component rational curve
that has a tacnode. Thus, $Z\cap W_{ij}\!=\!\eset$.
\end{lmm}

\noindent
{\it Proof:}
(1) The family $\tilde{\cal F}$ corresponding to the case of~$W_{41}$
is obtained by blowing up the family ${\cal F}$
of Lemma~\ref{nc_lmm1}
at two smooth points, $p_1$ and $p_2$, of~$D_1\!\subset\! C_{13}$.
On a neighborhood of $D_i^*\!\subset\!C_{41}$, we use local coordinates
$(t,z_i)\!\lra\!\big(t,p_i\!+\!tz_i,[1,z_i]\big)$.
Then,
\begin{gather*}
\psi_t\big|_{z_i}=a_it\big(1+o(1_t)\big)dz_i;\\
\mu_t^{-1}(L_i)=\big\{ z_{1,1}^{(i)}(t),\ldots,z_{1,d_1}^{(i)}(t),
z_{2,d_1+1}^{(i)}(t),\ldots,z_{2,d}^{(i)}(t)\big\}\subset C_t,
\quad
z_{\io,j}^{(i)}(t)=z_{\io,j}^{(i)}(0)+o(1_t);\\
\sum_{j=1}^{j=d_1}\int_{z^{(1)}_{1,j}(t)}^{z^{(2)}_{1,j}(t)}\psi_t+
\sum_{j=d_1+1}^{j=d}\int_{z^{(1)}_{2,j}(t)}^{z^{(2)}_{2,j}(t)}\psi_t=0
\qquad\forall t\!\in\!\De^*.
\end{gather*}
for some $a_1,a_2\!\in\!\Bbb{C}^*$, which depend on $D$,
but not on~$\mu|D_i$.
Proceeding as before, we obtain
\begin{gather}
a_1t\sum_{j=1}^{j=d_1}\big(z^{(1)}_{1,j}(0)-z^{(2)}_{1,j}(0)+o(1_t)\big)+
a_2t\sum_{j=d_1+1}^{j=d}\big(z^{(2)}_{2,j}(0)-z^{(1)}_{2,j}(0)+o(1_t)\big)=0
\qquad\forall t\!\in\!\De^*\Lra\notag\\
\label{nc_lmm2_e3}
a_1\sum_{j=1}^{j=d_1}z^{(1)}_{1,j}(0)+
a_2\sum_{j=d_1+1}^{j=d}z^{(1)}_{2,j}(0)=
a_1\sum_{j=1}^{j=d_1}z^{(2)}_{1,j}(0)+
a_2\sum_{j=d_1+1}^{j=d}z^{(2)}_{2,j}(0).
\end{gather}
Let $p^{(1)}_i$ and $p^{(2)}_i$ be the homogeneous polynomials
corresponding to $\mu|D_1$ and $\mu|D_2$, respectively.
Since~\e_ref{nc_lmm2_e3} holds for a dense subset of lines,
there exist $K\!=\!K(\mu)\!\in\!\Bbb{C}$ such that
\begin{equation}\label{nc_lmm2_e5}
a_1\frac{c_0p_{0,d_1-1}^{(1)}+c_1p_{1,d_1-1}^{(1)}+c_2p_{2,d_1-1}^{(1)}}
{c_0p_{0,d_1}^{(1)}+c_1p_{1,d_1}^{(1)}+c_2p_{2,d_1}^{(1)}}
+a_2\frac{c_0p_{0,d_2-1}^{(2)}+c_1p_{1,d_2-1}^{(2)}+c_2p_{2,d_2-1}^{(2)}}
{c_0p_{0,d_2}^{(2)}+c_1p_{1,d_2}^{(2)}+c_2p_{2,d_2}^{(2)}}
=K,
\end{equation}
for all $(c_0,c_1,c_2)\!\in\!\Bbb{C}^3\!-\!\{0\}$.
Since $\mu$ maps the nodes of $D_1$ and $D_2$
to the same point,
$$\big(p_{0,d_1}^{(1)},p_{1,d_1}^{(1)},p_{2,d_1}^{(1)}\big)
=\ka\big(p_{0,d_2}^{(2)},p_{1,d_2}^{(2)},p_{2,d_2}^{(2)}\big)$$
for some $\ka\!\in\!\Bbb{C}$.
Using this equation, it is easy to see that 
condition~\e_ref{nc_lmm2_e5} is equivalent to saying that
$\mu$ maps the singular points of $D_1$ and $D_2$
into a tacnode of its image.
Thus, the image of very element of $\ov{W}_2(d)\cap W_{41}$
is a two-component curve with a tacnode.\\
(2) 
Nearly the same argument applies to $W_{51}$.
In this case, an extra blowup is required,
and we will have $a_1\!=\!a_2\!=\!a$.

\begin{lmm}
\label{nc_lmm5}
The image of every element of $\ov{W}_2(d)\cap W_{53}$ 
is a two-component rational curve such that
both components have a cusp at one of the nodes of 
the image curve. Thus, $Z\cap W_{53}\!=\!\eset$.
\end{lmm}

\noindent
{\it Proof:} 
The proof is a minor modification of the proof of Lemma~\ref{nc_lmm1}.
The central fiber of ${\cal F}$ in this case is~$C_{32}$.
We can then choose $\psi\!\in\! H^0\big(C_{32};\om_{C_{32}}\big)$
such that the restriction of $\psi$ to the right vertical component
(in the diagram) is zero.
In terms of coordinates $(t,w_1)$ and $(t,w_2)$
near the smooth points $p_1$ and $p_2$
of the two vertical components, we will have
$$\psi_t\big|_{w_1}=a\big(1+o(1_{(t,w)})\big)dw_1
\quad\hbox{and}\quad
\psi_t\big|_{w_2}=o(1_t)dw_2,$$
for some $a\!\in\!\Bbb{C}^*$.
Proceeding as above, we then conclude that 
$\mu$ maps the node of $D_1\!\subset\! C_{53}$ to a cusp of~$\mu(D_1)$.
The same argument applies to~$\mu|D_2$.

\subsection{The Remaining Cases}

\noindent
The arguments in the previous subsection look very much
like the arguments in~\cite{P} and~\cite{KQR}.
However, some differences appear in this subsection.

\begin{lmm}
\label{contr_lmm1}
If $\big[\mu\!:(D,p_1,\ldots,p_{3d-2})\big]\!\in\! \ov{W}_2(d)\cap W_{13}$,
the image of $\mu$ is a tacnodal rational curve
and $\mu$ maps the nodes of $D$ to a tacnode of $\mu(D)$.
\end{lmm}

\noindent
{\it Proof:}
(1) We use coordinates $(t,w)$ near $D_1^*\!\subset\! C_{13}$
such that the two nodes of $D_1$ correspond to $w\!=\!0$
and $w\!=\!\i$.
Let $\psi_t\!\in\! H^0\big(C_t;\om_{C_t}\big)$ be such that
$$\psi_t\big|_w=\big(1+o(1_t)\big)\frac{dw}{w}.$$
Proceeding as above, we obtain
\begin{gather}
\mu_t^{-1}(L_i)=\big\{ w_1^{(i)}(t),\ldots,w_d^{(i)}(t)\big\}\subset C_t,
\quad w_j^{(i)}(t)=w_j^{(i)}(0)+o(1_t);\notag\\
\sum_{j=1}^{j=d}\int_{w^{(1)}_j(t)}^{w^{(2)}_j(t)}\psi_t
=0\in \Bbb{C}\big/2\pi i\Bbb{Z}
\qquad\forall t\!\in\!\De^*;\notag\\
\label{contr_lmm1_e5a}
\prod_{j=1}^{j=d}w_j^{(1)}(0)=\prod_{j=1}^{j=d}w_j^{(2)}(0)\equiv K;\\
\label{contr_lmm1_e5b}
\big(p_{0,0},p_{1,0},p_{2,0}\big)=
K\big(p_{0,d},p_{1,d},p_{2,d}\big).
\end{gather} 
for some $K\!=\!K(\mu)\!\in\!\Bbb{C}$.
Condition~\e_ref{contr_lmm1_e5b} on the coefficients
of the homogeneous polynomials corresponding to $\mu|D_1$
follows from the fact that \e_ref{contr_lmm1_e5a}
holds for a dense subset of lines in~$\PP$.
However, \e_ref{contr_lmm1_e5b} by itself tells us nothing new
about $\mu|D_1$, since we already know that $\mu$
maps the nodes of $D_1$ to the same point.\\
(2)  We instead consider the limit of 
the left-hand side of \e_ref{contr_lmm1_e5a}
as $L_1$ approaches the line tangent to 
the branch $w\!=\!0$ of~$\mu(D)$.
If the node $\mu(0)$ of $\mu(D)$ is simple,
two of the numbers $w_j^{(1)}(0)$ tend to $0$
and one to $\i$, all at comparable rates.
Thus, we must have $K\!=\!0$.
By the same argument, $K\!=\!\i$.
This means
$$p_{0,0}=p_{1,0}=p_{2,0}=p_{0,d}=p_{1,d}=p_{2,d}=0.$$
If $[u,v]$ are homogeneous coordinates on $E^{(1)}$
with $w=\frac{v}{u}$, it follows that $uv$ divides
all three homogeneous polynomials $p_0,p_1,p_2$,
i.e.~$\mu|D$ has degree at most $d\!-\!2$, not~$d$,
contrary to the assumption.
Thus, $\mu(0)\!=\!\mu(\i)$ has to be a tacnode of~$\mu(D)$
if \hbox{$\big[\mu\!:(D,p_1,\ldots,p_{3d-2})\big]\!\in\! 
\ov{W}_2(d)\cap W_{13}$}.

\begin{lmm}
\label{nc_lmm6}
The image of every element 
$\big[\mu\!:(D,p_1,\ldots,p_{3d-2})\big]\!\in\!\ov{W}_2(d)\cap W_{32}$
has a tacnode at $\mu(p_i)$ for some $i\!=\!1,\ldots,3d\!-\!2$.
If $\big[\mu\!:(D,p_1,\ldots,p_{3d-2})\big]\!\in\!\ov{W}_2(d)\cap W_{24}$
and the degree of $\mu|D_1$ is less than~$d$,
then $\mu(D)$ is a two-component rational tacnodal curve.
Thus, $Z\cap W_{ij}\!=\!\eset$ in both cases.
\end{lmm}

\noindent
{\it Proof:} Since the proof of Lemma~\ref{contr_lmm1}
carries over to the case of $W_{32}$ with no change,
the first claim is clear.
For the second claim, we use coordinates $(t,w)$ and
$(t,z)$ as in the proofs of Lemmas~\ref{nc_lmm1}
and~\ref{contr_lmm1}.
Then,
\begin{gather*}
\psi_t\big|_w=\big(1+o(1_t)\big)\frac{dw}{w},\qquad
\psi_t\big|_z=o(1_t);\\
\mu_t^{-1}(L_i)=
\big\{ z_1^{(i)}(t),\ldots,z_{d_1}^{(i)}(t),
w_{d_1+1}^{(i)}(t),\ldots,w_d^{(i)}(t)\big\}\subset C_t,
\qquad\hbox{with}\\ 
z_j^{(i)}(t)=z_j^{(i)}(0)+o(1_t),\quad
w_j^{(i)}(t)=w_j^{(i)}(0)+o(1_t);\\
\sum_{j=1}^{j=d_1}\int_{z^{(1)}_j(t)}^{z^{(2)}_j(t)}\psi_t
+\sum_{j=d_1+1}^{j=d}\int_{w^{(1)}_j(t)}^{w^{(2)}_j(t)}\psi_t
=0\in \Bbb{C}\big/2\pi i\Bbb{Z}
\qquad\forall t\!\in\!\De^*;\\
\prod_{j=d_1+1}^{j=d}w_j^{(1)}(0)=
\prod_{j=d_1+1}^{j=d}w_j^{(2)}(0).
\end{gather*} 
The last identity implies that $\mu|D_2$ has a tacnode.
The remaining claim of the lemma follows from the first two
and Corollary~\ref{nc_lmm1c}.

\begin{lmm}
\label{nc_lmm7}
If $\big[\mu\!:(D,p_1,\ldots,p_{3d-2})\big]\!\in\! \ov{W}_2(d)\cap W_{31}$
and the degree of $\mu|D_1$ is~$d$,
$\mu(D)$ has a tacnode at $\mu(p_i)$ for some $i\!=\!1,\ldots,3d\!-\!2$.
If the degree of $\mu|D_1$ is less than~$d$,
$\mu(D)$ is a two-component tacnodal rational curve.
Thus, $Z\cap W_{31}\!=\!\eset$.
\end{lmm}

\noindent
{\it Proof:}
The proof of Lemma~\ref{contr_lmm1} applies to the first case 
with no change.
For the second case, we use coordinate $(t,w_1)\!=\!(t,w)$ and
$(t,w_2)$ analogous to~$(t,w)$,
such that $w_1\!=\!\i$ and $w_2\!=\!\i$ are identified in~$C_{31}$.
Since the residues of 
$\psi\!\in\! H^0\big(\tilde{C}_0;\om_{\tilde{C}_0}\big)$
at $w_1\!=\!\i$ and $w_2\!=\!\i$
add up to zero, $\psi|D_2\!=\!-\frac{dw_2}{w_2}$.
Thus, proceeding as in the proof of Lemma~\ref{contr_lmm1},
we obtain
\begin{gather}
\prod_{j=1}^{j=d_1}w_{1,j}^{(1)}(0)
\cdot \Big(\prod_{j=d_1+1}^{j=d}w_{2,j}^{(1)}(0)\Big)^{-1}
=\prod_{j=1}^{j=d_1}w_{1,j}^{(2)}(0)
\cdot \Big(\prod_{j=d_1+1}^{j=d}w_{2,j}^{(2)}(0)\Big)^{-1}\equiv K;\notag\\
\label{nc_lmm4_e3}
\frac{c_0p_{0,0}^{(1)}+c_1p_{1,0}^{(1)}+c_2p_{2,0}^{(1)}}
{c_0p_{0,d_1}^{(1)}+c_1p_{1,d_1}^{(1)}+c_2p_{2,d_1}^{(1)}}
\cdot
\frac{c_0p_{0,d_2}^{(2)}+c_1p_{1,d_2}^{(2)}+c_2p_{2,d_2}^{(2)}}
{c_0p_{0,0}^{(2)}+c_1p_{1,0}^{(2)}+c_2p_{2,0}^{(2)}}
=K\quad\forall(c_0,c_1,c_2)\!\in\!\Bbb{C}^3\!-\!\{0\},
\end{gather}
for some $K\!\in\!\Bbb{C}$.
Since $\mu\big(w_2\!=\!\i\big)=\mu\big(w_1\!=\!\i\big)$,
$$\big(p_{0,d_1}^{(1)},p_{1,d_1}^{(1)},p_{2,d_1}^{(1)}\big)
=\ka\big(p_{0,d_2}^{(2)},p_{1,d_2}^{(2)},p_{2,d_2}^{(2)}\big)$$
for some $\ka\!\in\!\Bbb{C}^*$.
Thus, as a condition on $\mu$, \e_ref{nc_lmm4_e3} is equivalent to
$$\big(p_{0,0}^{(1)},p_{1,0}^{(1)},p_{2,0}^{(1)}\big)
=K\big(p_{0,0}^{(2)},p_{1,0}^{(2)},p_{2,0}^{(2)}\big)$$
for some $K\!\in\!\Bbb{C}$.
Suppose $\mu\big(w_2\!=\!\i\big)=\mu\big(w_1\!=\!0\big)$ 
is not a tacnode of~$\mu(D)$.
Then as in (2) of the proof of Lemma~\ref{contr_lmm1},
we conclude that
$$p_{0,0}^{(1)}=p_{1,0}^{(1)}=p_{2,0}^{(1)}\
=p_{0,0}^{(2)}=p_{1,0}^{(2)}=p_{2,0}^{(2)}.$$
This means $\mu|D_1$ and $\mu|D_2$
have degrees at most $d_1\!-\!1$ and $d_2\!-\!1$, respectively,
contrary to the assumption.

\section{Proof of Proposition~\ref{extra_contr}}
\label{extra_contr_sec}

\noindent
By Lemma~\ref{contr_lmm1}, if
$\big[\mu\!:(D,p_1,\ldots,p_{3d-2})\big]\!\in\!Z\cap W_{13}$,
$\mu$ maps the nodes of $D$ into the tacnode of~$\mu(D)$.
We now prove the converse and determine the multiplicity
with which the number $T_d$ enters into $[\pi^{-1}(C_0)]\cdot Z$.

\begin{lmm}
\label{contr_lmm2}
Suppose $C_0'$ is a tacnodal rational curve and 
$\eta\!:{\cal W}\!\lra\!{\cal B}$
is a local deformation space for~$C_0$.
Let $q_1,\ldots,q_{3d-2}$ be points in general position 
in~$\PP$ and $f\!:C_0'\!\lra\!\PP$ be a map of degree~$d$ 
passing through the $(3d\!-\!2)$ points.
Then there exists a map $\tilde{f}\!:{\cal W}\!\lra\!\PP$,
perhaps after shrinking~${\cal B}$, such that 
$\tilde{f}|C_0'\!=\!f$ and 
$\tilde{f}\big|\eta^{-1}(t)$ passes through the $(3d\!-\!2)$ points.
\end{lmm}

\noindent
{\it Proof:}
Since $T_d\!=\!0$ for $d\!\le\!3$, we can assume $d\!\ge\!3$.
Then $H^1\big(C_0';f^*{\cal O}_{\PP}(1)\big)\!=\!0$,
or equivalently $H^1\big(C_0';f^*T\PP\big)\!=\!0$.
Thus, there is no obstruction to extending $f$ to a neighborhood
of $C_0'$ in~${\cal W}$.

\begin{crl}
\label{contr_lmm2c}
Suppose $\big[\mu\!:(D,p_1,\ldots,p_{3d-2})\big]\!\in\! W_{13}$,
$\mu(p_i)\!=\!q_i$ for all $i\!=\!1,\ldots,3d\!-\!2$, and
$\mu$ maps the nodes of $D_1$ to the tacnode of~$\mu(D)$.
Then $\big[\mu\!:(D,p_1,\ldots,p_{3d-2})\big]\!\in\!\ov{W}_2(d)$.
\end{crl}

\noindent
{\it Proof:}
We apply Lemma~\ref{contr_lmm2} to the normalization
$f\!:C_0'\!\lra\!\mu(D)$ of $\mu(D)$ at the simple nodes.
Let $C_t$ be a family of a rational curves identified at
two pairs of points, i.e.

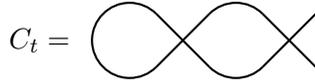
\begin{pspicture}(0,-1)(10,1)
\psarc(6,0){.5}{45}{315} 
\psline(6.354,-.354)(7.061,.354)\psline(6.354,.354)(7.061,-.354)
\psarc(7.414,0){.5}{45}{135}\psarc(7.414,0){.5}{225}{315}
\psline(7.768,-.354)(8.475,.354)\psline(7.768,.354)(8.475,-.354)
\rput(4.8,0){$C_t=$}
\end{pspicture}

\noindent
As the nodes of $C_t$ come together, $C_t$ approaches $C_0'$ in ${\cal B}$.
For all $t\!\neq\!0$ sufficiently small, 
let $f_t\!:C_t\!\lra\!\PP$ be the maps provided by Lemma~\ref{contr_lmm2}.
Then $f_t(C_t)$ converges to~$f(C_0')$.
Furthermore, $C_t$ converges to $C_0$ in $\ov{\frak M}_2$.
Thus, if 
$$\lim_{t\lra0}\big[f_t\!:(C_t,f_t^{-1}(q_1),\ldots,f_t^{-1}(q_{3d-2}))\big]
=\big[\mu'\!:(D',p_1',\ldots,p_{3d-2}')\big]
\in\ov{\frak M}_2(d),$$
$D'$ must be one of the curves $C_{ij}$ of Figure~1,
and $\mu'(D')$ is a tacnodal rational curve.
By Propositions~\ref{no_contr1} and~\ref{no_contr2},
we conclude that 
$$\big[\mu\!:(D,p_1,\ldots,p_{3d-2})\big]
=\big[\mu':(D',p_1',\ldots,p_{3d-2}')\big]
\in\ov{\frak M}_2(d).$$

\begin{lmm}
\label{contr_l3}
The contribution of $W_{13}$ to $[\pi^{-1}(C_0)]\cdot Z$
is $6T_d$.
\end{lmm}

\noindent
{\it Proof:} Suppose 
$\big[\mu\!:(D,p_1,\ldots,p_{3d-2})\big]\!\in\!Z\cap W_{13}$.
Given a fixed complex structure $j$ on $\Si$ such that
$(\Si,j)$ is very close to $[C_0]$ in $\ov{\frak M}_2$,
we need to determine the number maps 
$\mu_j\!:\Si\!\lra\!\PP$ close to~$\mu$.
By Corollary~\ref{contr_lmm2c}, there exists a family of
curves 
$\tilde{\eta}\!:\tilde{F}\!\lra\!\De$ and 
of maps $\tilde{\mu}\!:\tilde{F}\!\lra\!\PP$
restricting to $\mu$ on the central fiber~$D$.
There are six automorphisms of $C_0$ that preserve  its components.
Corresponding to these automorphisms and $(\tilde{F},\tilde{\eta})$,
we obtain six maps $\mu_j\!:\Si\!\lra\!\PP$.
None of these maps are equivalent, since we did not switch 
the two components of~$C_0$.

\end{document}